\documentclass[10pt]{article}
\usepackage[mathscr]{eucal}
\usepackage{mathrsfs}
\usepackage{amssymb}
\usepackage{calligra}
\usepackage{fontenc}
\usepackage{amsbsy} 
\usepackage{graphicx}
\graphicspath{%
    {converted_graphics/}
    {/}
}
\begin{document}

\newcommand{\ts}{{\tilde{\sf s}}}
\newcommand{\sfv}{{\sf v}}
\newcommand{\sfw}{{\sf w}}
\newcommand{\simge}{\ba{cc}\vspace*{-2.4mm}>\\ \sim\ea }
\newcommand{\simle}{\ba{cc}\vspace*{-2.4mm}<\\ \sim\ea }
\newcommand{\Cdot}{\!\cdot\!}
\newcommand{\sq}{{$\sqcap\!\!\!\!\sqcup$}}
\newcommand{\Eu}{{\rm I\,\!\! E}}
\newcommand{\Io}{\Int{\Omega}{}}
\newcommand{\Id}{\Int{\cald}{}}
\newcommand{\Div}{\mbox{\rm div}\,}
\newcommand{\tr}{\mbox{\rm tr}\,}
\newcommand{\grad}{\mbox{\rm grad}\,}
\newcommand{\supp}{\mbox{\rm supp}\,}
\newcommand{\curl}{\mbox{\rm curl}\,}
\newcommand{\Ido}{\Int{\partial\Omega}{}}
\newcommand{\IdS}{\Int{\Sigma}{}}
\newcommand{\Oint}[2]{{\displaystyle \oint_{#1}^{#2}}}
\newcommand{\Int}[2]{{\displaystyle \int_{ #1}^{ #2}}}
\newcommand{\Lim}[1]{{\displaystyle \lim_{ #1}}}
\newcommand{\Limsup}[1]{{\displaystyle \limsup_{\footnotesize #1}}}
\newcommand{\Liminf}[1]{{\displaystyle \liminf_{\footnotesize #1}}}
\newcommand{\Sup}[1]{{\displaystyle \sup_{#1}}}
\newcommand{\Inf}[1]{{\displaystyle \inf_{#1}}}
\newcommand{\Max}[1]{{\displaystyle \max_{#1}}}
\newcommand{\Min}[1]{{\displaystyle \min_{#1}}}
\newcommand{\Sum}[2]{{\displaystyle \sum_{#1}^{#2}}}
\newcommand{\Prod}[2]{{\displaystyle \prod_{#1}^{#2}}}
\newcommand{\BCup}[2]{{\displaystyle \bigcup_{#1}^{#2}}}
\newcommand{\BCap}[2]{{\displaystyle \bigcap_{#1}^{#2}}}
\newcommand{\Frac}[2]{\displaystyle{\frac{\displaystyle{#1}}{\displaystyle{#2}}}}
\newcommand{\norm}[1]{\left\|{#1}\right\|}
\newcommand{\Norm}[1]{\langle\langle{#1}\rangle\rangle_q}
\newcommand{\No}[1]{\langle\!\langle{#1}\rangle\!\rangle}
\newcommand{\NO}[1]{{\langle{#1}\rangle}_{\lambda,q}}
\newcommand{\beea}{\begin{eqnarray}}
\newcommand{\eeea}{\end{eqnarray}}
\newcommand{\ms}{\medskip\smallskip}
\newcommand{\bs}{\bigskip}
\newcommand{\ps}{\par\smallskip}
\newcommand{\bfe}{{\mbox{\boldmath $e$}} }
\newcommand{\pni}{{\par\noindent}}
\newcommand{\bfq}{{\mbox{\boldmath $q$}} }
\newcommand{\bfz}{{\mbox{\boldmath $z$}} }
\newcommand{\0}{{\mbox{\boldmath $0$}} }
\newcommand{\LE}{\!\!\!&\le&\!\!\!}
\newcommand{\BL}[1]{{\par\smallskip{\bf Lemma #1.}}}
\newcommand{\BT}[1]{{\par\smallskip{\bf Theorem #1.}}}
\newcommand{\Ln}{[\!|}
\newcommand{\Rn}{|\!]}
\newcommand{\n}[1]{{\Ln{#1}\Rn}} 
\newcommand{\nq}[1]{{\Ln{#1}\Rn}_{q}} 
\newcommand{\nqr}[1]{{\Ln{#1}\Rn}_{q,r}} 
\newcommand{\Nq}[1]{{\langle{#1}\rangle}_{q}} 
\newcommand{\Nql}[1]{{\langle{#1}\rangle}_{\lambda,q}} 
\newcommand{\Nqr}[1]{{\langle{#1}\rangle}_{q,r}}
\newcommand{\N}[1]{{|\!\!|\!\!|\,{#1}\,|\!|\!\!|_2}}
\newcommand{\EA}[2]{$$#1$$%
\vspace{-6.mm}
\begin{equation}
\end{equation}
\vspace{-6.mm}
$$
#2
\setlength{\belowdisplayskip}{3mm}
\setlength{\belowdisplayshortskip}{3mm}
$$
}
\newcommand{\A}[2]{$$#1$$%
\vspace{-4.mm}
$$
#2
\setlength{\belowdisplayskip}{3mm}
\setlength{\belowdisplayshortskip}{3mm}
$$
}
\newcommand{\BF}{\begin{footnotesize}}
\newcommand{\EF}{\end{footnotesize}}
\setlength{\jot}{.15in}
\newcommand{\pde}[2]{{\displaystyle \frac{\mbox{$\partial #1$}}{\mbox{$\partial #2$}}}}
\newcommand{\ode}[2]{{\displaystyle \frac{\mbox{$d #1$}}{\mbox{$d #2$}}}}
\newcommand{\f}[2]{\frac{\mbox{$#1$}}{\mbox{$ #2$}}}
\newcommand{\bi}{\begin{itemize}}
\newcommand{\ei}{\end{itemize}}
\newcommand{\ed}{\end{document}}
\newcommand{\be}{\begin{equation}}
\newcommand{\ba}{\begin{array}}
\newcommand{\ea}{\end{array}}
\newcommand{\ee}{\end{equation}}
\newcommand{\eeq}[1]{\label{eq:#1}\end{equation}}
\newcommand{\real}{{\mathbb R}}
\newcommand{\compl}{{\mathbb C}}
\def\Id{\mbox{\boldmath $1$}}
\def\zero{\mbox{\boldmath $0$}}
\newcommand{\PP}{{\rm I\!\!\,P}}
\newcommand{\nat}{{\mathbb N}}
\newcommand{\bfpsi}{\mbox{\boldmath $\psi$}}
\newcommand{\bfchi}{\mbox{\boldmath $\chi$}}
\newcommand{\bfomega}{\mbox{\boldmath $\omega$}}
\newcommand{\bfome}{\mbox{\boldmath $\varpi$}}
\newcommand{\bfvaromega}{\mbox{\boldmath $\varpi$}}
\newcommand{\bfOmega}{\mbox{\boldmath $\Omega$}}
\newcommand{\bfTheta}{\mbox{\boldmath $\Theta$}}
\newcommand{\bfxi}{\mbox{\boldmath $\xi$}}
\newcommand{\bfmu}{\mbox{\boldmath $\mu$}}
\newcommand{\bfx}{\mbox{\boldmath $x$}}
\newcommand{\bfy}{\mbox{\boldmath $y$}}
\newcommand{\bfPsi}{\mbox{\boldmath $\Psi$}}
\newcommand{\bfphi}{\mbox{\boldmath $\varphi$}}
\newcommand{\bfhi}{\mbox{\boldmath $\phi$}}
\newcommand{\bfPhi}{\mbox{\boldmath $\Phi$}}
\newcommand{\bfv}{{\mbox{\boldmath $v$}} }
\newcommand{\bfu}{{\mbox{\boldmath $u$}} }
\newcommand{\bfsf}{{\mbox{\footnotesize\boldmath $s$}} }
\newcommand{\bfuf}{{\mbox{\footnotesize\boldmath $u$}} }
\newcommand{\bfw}{{\mbox{\boldmath $w$}} }
\newcommand{\bff}{{\mbox{\boldmath $f$}} }
\newcommand{\bfa}{{\mbox{\boldmath $a$}} }
\newcommand{\bfi}{{\mbox{\boldmath $i$}} }
\newcommand{\bfj}{{\mbox{\boldmath $j$}} }
\newcommand{\bfc}{{\mbox{\boldmath $c$}} }
\newcommand{\bfo}{{\mbox{\boldmath $o$}} }
\newcommand{\bfp}{{\mbox{\boldmath $p$}} }
\newcommand{\bfkp}{{\mbox{\footnotesize{\boldmath $k$}}} }
\newcommand{\bfka}{{\mbox{\footnotesize{\boldmath $k^*$}}} }
\newcommand{\bft}{{\mbox{\boldmath $t$}} }
\newcommand{\bfd}{{\mbox{\boldmath $d$}} }
\newcommand{\bfl}{{\mbox{\boldmath $l$}} }
\newcommand{\bfr}{{\mbox{\boldmath $r$}} }
\newcommand{\bfk}{{\mbox{\boldmath $k$}} }
\newcommand{\bfA}{{\mbox{\boldmath $A$}} }
\newcommand{\bfS}{{\mbox{\boldmath $S$}} }
\newcommand{\bfO}{{\mbox{\boldmath $O$}} }
\newcommand{\bfM}{{\mbox{\boldmath $M$}} }
\newcommand{\bfP}{{\mbox{\boldmath $P$}} }
\newcommand{\bfB}{{\mbox{\boldmath $B$}} }
\newcommand{\bfR}{{\mbox{\boldmath $R$}} }
\newcommand{\bfC}{{\mbox{\boldmath $C$}} }
\newcommand{\bfD}{{\mbox{\boldmath $D$}} }
\newcommand{\bfQ}{{\mbox{\boldmath $Q$}} }
\newcommand{\bfZ}{{\mbox{\boldmath $Z$}} }
\newcommand{\bfG}{{\mbox{\boldmath $G$}} }
\newcommand{\bfE}{{\mbox{\boldmath $E$}} }
\newcommand{\bfX}{{\mbox{\boldmath $X$}} }
\newcommand{\bfY}{{\mbox{\boldmath $Y$}} }
\newcommand{\bfH}{{\mbox{\boldmath $H$}} }
\newcommand{\bfI}{{\mbox{\boldmath $I$}} }
\newcommand{\bfJ}{{\mbox{\boldmath $J$}} }
\newcommand{\bfN}{{\mbox{\boldmath $N$}} }
\newcommand{\bfh}{{\mbox{\boldmath $h$}} }
\newcommand{\bfm}{{\mbox{\boldmath $m$}} }
\newcommand{\bfone}{{\mbox{\boldmath $1$}} }
\newcommand{\hs}{{\rm I}\!\!\,{\rm R}^3_+}
\newcommand{\cala}{{\cal A}}
\newcommand{\calb}{{\cal B}}
\newcommand{\calc}{{\cal C}}
\newcommand{\cald}{{\cal D}}
\newcommand{\cale}{{\cal E}}
\newcommand{\calf}{{\cal F}}
\newcommand{\calg}{{\cal G}}
\newcommand{\calh}{{\cal H}}
\newcommand{\cali}{{\cal I}}
\newcommand{\calj}{{\cal J}}
\newcommand{\calk}{{\cal K}}
\newcommand{\call}{{\cal L}}
\newcommand{\calm}{{\cal M}}
\newcommand{\caln}{{\cal N}}
\newcommand{\calo}{{\cal O}}
\newcommand{\calp}{{\cal P}}
\newcommand{\calq}{{\cal Q}}
\newcommand{\calr}{{\cal R}}
\newcommand{\cals}{{\cal S}}
\newcommand{\calt}{{\cal T}}
\newcommand{\calu}{{\cal U}}
\newcommand{\calv}{{\cal V}}
\newcommand{\calx}{{\cal X}}
\newcommand{\caly}{{\cal Y}}
\newcommand{\calw}{{\cal W}}
\newcommand{\calz}{{\cal Z}}
\newcommand{\bfsigma}{\mbox{\boldmath $\sigma$}}
\newcommand{\bfSigma}{\mbox{\boldmath $\Sigma$}}
\newcommand{\bftau}{\mbox{\boldmath $\tau$}}
\newcommand{\bfeta}{\mbox{\boldmath $\eta$}}
\newcommand{\bfT}{{\mbox{\boldmath $T$}} }
\newcommand{\bfV}{{\mbox{\boldmath $V$}} }
\newcommand{\bfU}{{\mbox{\boldmath $U$}} }
\newcommand{\bfW}{{\mbox{\boldmath $W$}} }
\newcommand{\bfF}{{\mbox{\boldmath $F$}} }
\newcommand{\bfK}{{\mbox{\boldmath $K$}} }
\newcommand{\bfL}{{\mbox{\boldmath $L$}} }
\newcommand{\bfb}{{\mbox{\boldmath $b$}} }
\newcommand{\bfg}{{\mbox{\boldmath $g$}} }
\newcommand{\bfn}{{\mbox{\boldmath $n$}} }
\newcommand{\bfs}{{\mbox{\boldmath $s$}} }
\newcommand{\cf}{{\it cf.} }
\newcommand{\io}{\int_\Omega}
\newcommand{\1}{\item[({\it i})]}
\newcommand{\2}{\item[({\it ii})]}
\newcommand{\3}{\item[({\it iii})]}
\newcommand{\4}{\item[({\it iv})]}
\newcommand{\5}{\item[({\it v})]}
\newcommand{\6}{\item[({\it vi})]}
\newcommand{\7}{\item[({\it vii})]}
\newcommand{\8}{\item[({\it viii})]}
\newcommand{\9}{\item[({\it xi})]}
\newcommand{\ido}{\int_{\partial\Omega}}
\newcommand{\half}{\mbox{$\frac{1}{2}$}}
\def\parallel{\|}
\def\mid{|}
\def\Bbb R{\real}
\def\hat{\widehat}
\def\tilde{\widetilde}
\def\bar{\overline}
\newcommand{\threehalves}{3\over 2}
\newcommand{\bfPi}{\mbox{\boldmath $\Pi$}}
\newcommand{\bfXi}{\mbox{\boldmath $\Xi$}}
\newcommand{\bfalpha}{\mbox{\boldmath $\alpha$}}
\newcommand{\bfbeta}{\mbox{\boldmath $\beta$}}
\newcommand{\bfgamma}{\mbox{\boldmath $\gamma$}}
\newcommand{\bfdelta}{\mbox{\boldmath $\delta$}}
\newcommand{\bfzeta}{\mbox{\boldmath $\zeta$}}
\newcommand{\bfUpsilon}{\mbox{\boldmath $\Upsilon$}}
\newcommand{\bfGamma}{\mbox{\boldmath $\Gamma$}}
\newcommand{\bfcala}{\mbox{\boldmath ${\cal A}$}}
\newcommand{\bfcalm}{\mbox{\boldmath ${\cal M}$}}
\newcommand{\bfcaln}{\mbox{\boldmath ${\cal N}$}}
\newcommand{\bfcalq}{\mbox{\boldmath ${\cal Q}$}}
\newcommand{\bfcalb}{\mbox{\boldmath ${\cal B}$}}
\newcommand{\bfcalc}{\mbox{\boldmath ${\cal C}$}}
\newcommand{\bfcali}{\mbox{\boldmath ${\cal I}$}}
\newcommand{\bfcalg}{\mbox{\boldmath ${\cal G}$}}
\newcommand{\bfcalh}{\mbox{\boldmath ${\cal H}$}}
\newcommand{\bfcalk}{\mbox{\boldmath ${\cal K}$}}
\newcommand{\bfcalt}{\mbox{\boldmath ${\cal T}$}}
\newcommand{\bfcalx}{\mbox{\boldmath ${\cal X}$}}
\newcommand{\bfcall}{\mbox{\boldmath ${\cal L}$}}
\newcommand{\bfcalf}{\mbox{\boldmath ${\cal F}$}}
\newcommand{\bfcalr}{\mbox{\boldmath ${\cal R}$}}
\newcommand{\bfcals}{\mbox{\boldmath ${\cal S}$}}
\newcommand{\bfcalw}{\mbox{\boldmath ${\cal W}$}}
\newcommand{\bfcalu}{\mbox{\boldmath ${\cal U}$}}
\newcommand{\bfcalv}{\mbox{\boldmath ${\cal V}$}}
\newcommand{\bfcalz}{\mbox{\boldmath ${\cal Z}$}}
\pagenumbering{roman}
\newcommand{\art}[6]{{\I[{\sc #1,}] {#2}, {\it #3}, {\bf #4}, {#5} {[#6]}}}
\newcommand{\ED}{\end{description}}
\newcommand{\I}{\item }
\newcommand{\ra}{\rm a}
\newcommand{\rb}{\rm b}
\newcommand{\rc}{\rm c}
\newcommand{\Hsp}{{\rm I}\!\!\,{\rm R}^n_+}
\newcommand{\Hsn}{{\rm I}\!\!\,{\rm R}^n_-}
\newcommand{\po}[1]{\mbox{$\displaystyle \frac{\mbox{$\partial #1$}}
{\mbox{$\partial x_{1}$}}$}}
\newcommand{\PO}[1]{\mbox{$\displaystyle \frac{\mbox{$\partial #1$}}
{\mbox{$\partial y_{1}$}}$}}
\newcommand{\OP}{\left(\Delta+2\lambda\PO{}\right)}
\newcommand{\op}{\left(\Delta+2\lambda\po{}\right)}
\newcommand{\ft}[1]{
\Frac{1}{(2\pi)^{n/2}}\Int{{\Bbb R}^{n}}{}e^{i{\bf x}\cdot \bfxi}
#1(\xi)d\xi}
\newcommand{\Ft}[1]{
\Frac{1}{2\pi}\Int{{\Bbb R}^{2}}{}e^{i{x}\cdot \xi}
#1(\xi)d\xi}
\newcommand{\Z}{\item[({\it a})]}
\newcommand{\B}{\item[({\it b})]}
\newcommand{\C}{\item[({\it c})]}
\newcommand{\D}{\item[({\it d})]}
\newcommand{\E}{\item[({\it e})]}
\newcommand{\G}{\item[({\it g})]}
\newcommand{\Š}{\`e}
\newcommand{\…}{\`a}
\newcommand{\•}{\`o}
\newcommand{\—}{\`u}
\newcommand{\}{\`{\i}}
\def\tag{\renewcommand{\theequation}}
\newcommand{\Footnote}{~\footnote}
\newcommand{\ie}{{\it i.e.}}
\newcommand{\dist}{\mbox{\rm dist\,}}
\newcommand{\const}{\mbox{\rm const}}
\newcommand{\trace}{\mbox{\rm trace}}
\newcommand{\Bo}{\par\hfill{$\Box$}\par\noindent}
\newcommand{\Nor}[1]{\langle{#1}\rangle_q}
\newcommand{\vs}{\vspace*{.5cm}\par\noindent}
\newcommand{\Vs}{\vspace*{.6cm}\par\noindent}
\newcommand{\Vvs}{\vspace*{.7cm}\par\noindent}
\newcommand{\VVs}{\vspace*{.8cm}\par\noindent}
\newtheorem{definition}{Definition}[section]
\newcommand{\Bd}{\begin{definition}\begin{rm}}
\newcommand{\Ed}{\end{rm}\end{definition}}
\newtheorem{remark}{Remark}[section]
\newcommand{\Br}{\begin{remark}\begin{rm}}
\newcommand{\Er}{\end{rm}\end{remark}}
\newtheorem{proposition}{Proposition}[section]
\newcommand{\Bp}{\begin{proposition}\begin{sl}}
\newcommand{\EP}[1]{\end{sl}\label{proposition:#1}\end{proposition}}
\newcommand{\propref}[1]{{\rm Proposition \ref{proposition:#1}}}
\newcommand{\Bt}{\begin{theorem}\begin{sl}}
\newcommand{\Et}{\end{sl}\end{theorem}}
\newcommand{\Bl}{\begin{lemma}\begin{sl}}
\newcommand{\El}{\end{sl}\end{lemma}}
\newtheorem{theorem}{Theorem}[section]
\newtheorem{lemma}{Lemma}[section]
\newtheorem{corollary}{Corollary}[section]
\newcommand{\eqref}[1]{{\rm (\ref{eq:#1})}}
\newcommand{\Bc}{\begin{corollary}\begin{sl}}
\newcommand{\Ec}{\end{sl}\end{corollary}}
\newcommand{\ET}[1]{\end{sl}\label{theorem:#1}\end{theorem}}
\newcommand{\EDD}[1]{\end{rm}\label{definition:#1}\end{definition}}
\newcommand{\EL}[1]{\end{sl}\label{lemma:#1}\end{lemma}}
\newcommand{\theoref}[1]{{\rm Theorem \ref{theorem:#1}}}
\newcommand{\ER}[1]{\end{rm}\label{remark:#1}\end{remark}}
\newcommand{\EC}[1]{\end{sl}\label{corollary:#1}\end{corollary}}
\newcommand{\remref}[1]{{\rm Remark \ref{remark:#1}}}
\newcommand{\cororef}[1]{{\rm Corollary \ref{corollary:#1}}}
\newcommand{\lemmref}[1]{{\rm Lemma \ref{lemma:#1}}}
\newcommand{\essup}[1]{{\rm ess}\,{{\displaystyle \sup_{\hspace*{-5mm}{#1}}}}}

\renewcommand{\i}{{\rm i}}

\pagenumbering{arabic}
\newcommand{\QED}{{\par\hfill$\square$\par}}
\renewcommand{\thefootnote}{(\arabic{footnote})}
\title{On the Relation between Very Weak and Leray-Hopf Solutions\\ to Navier-Stokes Equations} 
\author{ Giovanni P. Galdi 
\thanks{Department of Mechanical Engineering and Materials Science, University of Pittsburgh, USA.
\,email: galdi@pitt.edu.
Work  partially supported by NSF DMS Grant-1614011.}}
\date{}

\maketitle
\begin{abstract} We prove a general result that implies that very weak solutions to the Cauchy problem for the Navier-Stokes equations must be, in fact, Leray-Hopf solutions if only their initial data are (solenoidal) with finite kinetic energy.    
 \end{abstract}
\renewcommand{\theequation}{\arabic{section}.\arabic{equation}}
\setcounter{section}{0}
\section{Introduction} 
%
%

We are concerned with the three-dimensional Cauchy problem\footnote{
We 
suppose, for simplicity, zero body force and,  without loss of generality, take the kinematic viscosity coefficient to be 1.} for the Navier-Stokes equations 
\be\ba{cc}\smallskip\left.\ba{ll}\smallskip
\partial_tv+v\cdot\nabla v=\Delta v-\nabla p\\
\Div v=0\ea\right\}\ \ \mbox{in}\,\ \real^3\times (0,\infty)
\\
v(x,0)=v_0(x)\ \ x\in\real^3\,,
\ea
\eeq{1.1}
where $v:\real^3\times  
[0,\infty)\to v(x,t)\in\real^3$ is the flow velocity field and $p$ the associated pressure field. \par
We recall that, for a given $v_0\in L^2_\sigma(\real^3)$, a corresponding {\em Leray--Hopf solution}  is a function $v$ with the property\footnote{$L^q_\sigma(\real^3)$ is the subspace of the Lebesgue space $L^q(\real^3)$ of divergence-free vector functions, and $C_w$ denotes the class of weakly continuous functions. Other notations are standard, like $H^{m,q}$,  for Sobolev spaces, with corresponding norm $\|\cdot\|_{m,q}$, $L^r(I;X)$, $I$  real interval, $X$ Banach space, for Bochner spaces, etc. }
$$
v\in C_w([0,T];L^2_\sigma(\real^3))\cap L^2(0,T;H^{1}(\real^3))\ \ \ \mbox{all $T>0$}\,,
$$
that solves \eqref{1.1} in a distributional sense, 
and satisfies the ``energy inequality:"
\be
\|v(t)\|_2^2+2\int_0^t\|\nabla v(\tau)\|_2^2\le\|v_0\|_2^2\,,\ \ \mbox{all $t\ge 0$}\,,
\eeq{1.3}
where $\|\cdot\|_q$ denotes the     $L^q(\real^3)$-norm. The existence of such a solution for an arbitrarily prescribed $v_0\in L^2_\sigma$  is known since the seminal papers of  Leray \cite{Leray} and Hopf \cite{Hopf}. 
\par 
Alongside with these solutions, there are the so-called {\em very weak} (or {\em mild}) {\em  solutions},  introduced, basically, in the pioneering work of Foias \cite{Foias}. Their properties were first consistently investigated  by Fabes {\em et al.} \cite{FJR}, and, more recently, by a number of authors, especially over the past two decades;  see, e.g., \cite{Kato,Amann,LiMa,FL,FarGaSoh,BeGa,Kuc,LR,FarKoSoh,FaR} and references therein. Precisely, a very weak solution is a field $v$ such that
\be
v\in L^r(0,T;L^s(\real^3))\,,\ \ \frac2r+\frac3s=1,\, \, s\in(3,\infty)\ \ \mbox{or}\ \ v\in C([0,T];L^3(\real^3)),\ \ T>0\,,
\eeq{h}
satisfying \eqref{1.1} in the following sense 
\be\ba{ll}\medskip
\Int{0}{T}\Int{\real^3}{}v\cdot(\partial_t\varphi+\Delta\varphi+v\cdot\nabla\varphi)=-\Int{\real^3}{}v_0\cdot\varphi(0)\\
\Int{0}{T}\Int{\real^3}{}v\cdot\nabla \phi=0\,,\ea
\eeq{1.5}
for arbitrary $\phi\in C_0^\infty(\real^3\times (0,T))$ and $\varphi\in\mathcal D_T:=\{\varphi\in C_0^\infty(\real^3\times [0,T)):\ \Div\varphi=0\}
$.
Notice that, at the outset, very weak solutions do not possess any  locally integrable derivative, and as a result, unlike Leray-Hopf's, may have infinite kinetic energy and overall dissipation. However,   under appropriate functional hypotheses on the initial data $v_0$, they exist, and are unique and  smooth at least in some time interval $[0,T^\star)$, with  $T^\star=\infty$ if the ``size" of $v_0$ is suitably restricted \cite{Kato,Amann,FarKoSoh}.
\par
On the other hand,  a  classical result in the Navier-Stokes theory  states that if a Leray-Hopf solution meets the requirement \eqref{h} (the so called {\em Prodi-Serrin-Ladyzhenskaya conditions}), it is then unique (in its class), smooth and obeys the {\em energy equality}, namely, \eqref{1.3} with the equality sign; e.g. \cite{Ga}. Nevertheless, as is well known, the existence of a Leray-Hopf solution satisfying \eqref{h} is an outstanding open question.  
\par 
Also motivated by the above considerations, several authors  have investigated under which assumptions --besides the necessary condition  $v_0\in L^2_\sigma$-- a very weak solution lies in the Leray-Hopf class.  The first contribution goes back to Foias \cite[Th\'eor\`eme 2]{Foias}, where the property is shown in case $v$ satisfies \eqref{h}$_1$ but with a {\em strict} inequality sign and, in addition, $v\in L^2(0,T;H^1)$. Fabes {\em et al.} \cite[Theorem (5.3)]{FJR} assume  \eqref{h}$_1$ with the supplemental condition $v_0\in L^s_\sigma$. The question was successively reconsidered by Kato 
 \cite{Kato} and, later on,  Giga \cite{Gi}.   In particular, in \cite[Proposition 1]{Gi} it is shown that if $v_0\in L^3_\sigma$, every corresponding very weak solution  $v$ satisfying {\em both} conditions in \eqref{h}, along with the further request  $t^{\frac1r}v\in BC([0,T);L^s(\real^3))$ and  $t^{\frac1r}\|v(t)\|_s\to 0$ as $t\to 0$, must be also Leray-Hopf. Very recently, Farwig and Riechwald \cite[Theorem 1]{FaR} have reached the same conclusion,  provided $v_0$ is in an appropriate negative Sobolev space and $v$ obeys \eqref{h}$_1$ along with the additional regularity property $v\in L^4(0,T;L^4(\real^3))$.
\par
Objective of this note is to show a general result for solutions to \eqref{1.5}, which  {\em implies, in particular, that every very weak solution is, in fact, Leray-Hopf provided only} $v_0\in L^2_\sigma$. Thus, a very weak solution will possess finite energy and overall dissipation in the time interval where it exists if only its initial energy is finite. \par More precisely, we shall prove the following.     
\Bt Let $v\in L^2_{\rm loc}(\real^3\times[0,T))$ satisfy
\eqref{1.5}
for some $v_0\in L^2_\sigma(\real^3)$ and all $\varphi\in \mathcal D_T$, $\phi\in C_0^\infty(\real^3\times (0,T))$.
Suppose that for all  small $\delta>0$, $v$ meets one of the following assumptions
\be\ba{ll}\smallskip
v\in L^{r}(\delta,T;L^s(\real^3))\,,\ \ r:=\frac{2s}{s-3}\,,\ \ s\in (3,\infty)\,,\\
v\in C([\delta,T];L^3(\real^3))\,,
\ea
\eeq{H}
and that $v\to v_0$ weakly in $L^2_\sigma(\real^3)$.
Then
$v$ is necessarily a Leray-Hopf solution. In addition, $ 
\lim_{t\to0^+}\|v(t)-v_0\|_2=0\,.
$
\ET{1}
\par
From this theorem  we can deduce a number of relevant  consequences.
\Bc {\rm (Energy Equality)} Let $v$ satisfy the assumption of \theoref{1}. Then $v$ obeys the energy equality (i.e. \eqref{1.3} with the equality sign) in the interval $[0,T]$.
\EC{1.1}
\Bc {\rm (Liouville-Type)} Let $v$ satisfy the assumption of \theoref{1} with $v_0\equiv 0$. Then $v\equiv 0$ in $[0,T]$.
\EC{1.1_0}
\Bc {\rm (Uniqueness)} Let $v$ obey  \eqref{h}.  Then $v$ is the only Leray--Hopf solution corresponding to $v_0$.
\EC{U}
\Bc {\rm (Regularity)} Let $v$ satisfy the assumptions of \theoref{1}. Then $v\in C^\infty(\real^3\times (0,T])$.
\EC{1.2}
The proof of \theoref{1}, given in Section 3, is based on ideas similar to those employed in \cite{GaPr} and is quite straightforward. In fact, it relies upon a simple duality argument, and a regularity result for solutions to a suitable linearization of the Navier-Stokes equations, presented in Section 2 (see \lemmref{2.4}). 
\section{Preparatory Results}\setcounter{equation}{0} 
Given a locally integrable $f:\real^3\times (0,T)\mapsto\real^3$ we define 
$$
f_{(\eta)}(x,t):=\Int0Tj_\eta(t-s)\big[\Int{\real^3}{}k_\eta(x-y)f(y,s)dy]
ds
\,,
$$
where 
$$
j_\eta(\tau):=\eta^{-1}j(\tau/\eta)\,,\ \ k_\eta(\xi):=\eta^{-1}k(\xi/\eta)\,,\ \ (\tau,\xi)\in \real\times\real^3\,, 
$$
with $j\in C_0^\infty(-1,1)$, and $k\in C_0^\infty(\real^3)$. 
We employ the usual notation
$
(f,g):=\int_{\real^3}f\cdot g\,,
$
$L^{p,q}:=L^p(0,T;L^q(\real^3))$,
and denote by $q':=q/(q-1)$ the H\"older conjugate of $q\in [1,\infty]$. 
For $p,q\in [1,\infty)$ we define the Banach space
$$
\mathscr W^{p,q}:=\big\{u\in L_{\rm 1oc}^1(\real^3\times (0,T)): u\in H^{1,p}(0,T;L^q_\sigma(\real^3))\cap L^p(0,T;H^{2,q}(\real^3))\big\}
$$
with corresponding norm
$$
\left\|u\right\|_{\mathscr W^{p,q}}:=\Big(\int_0^T\left(\|\partial_t u(t)\|_q^p+\|u(t)\|_{2,q}^p\right)\Big)^{\frac1p}\,.
$$
In the case $q=p$, we shall set $\mathscr W^{q,q}\equiv \mathscr W^{q}$. Finally, we put 
$$
\mathscr W^{p,q}_{0}:= \left\{\psi\in \mathscr W^{p,q}\,,\ \psi(T)=0\right\},$$ 
\smallskip\par
\Bl For any $p_i,q_i\in [1,\infty)$, $i=1,2$, the space $\mathcal D_T$ is dense in $\mathscr W^{p_1,q_1}_0\cap\mathscr W^{p_2,q_2}_0$. 
\EL{N1}
{\em Proof.} The proof can be achieved by an argument entirely analogous  to that of \cite[Lemma A.1]{GaPr}, and therefore it will be omitted.\par\hfill$\square$ 
\Bl Let  $\alpha,w_1$ be a given pair of functions satisfying  one of the following assumptions 
$$
{\rm (i)}\,\ \alpha,w_1\in L^{r,s}\,,\ r:=\frac{2s}{s-3}\,,\ s\in(3,\infty)\,;\ \ {\rm (ii)}\,\ \alpha,w_1\in L^{\infty,3}\,,
$$
and let $w_2\in L^{\infty,2}$.
 Consider the linear forms
$$
\mathscr T_i:\psi\in \mathscr W^{p,q}\mapsto \int_0^T(\alpha\cdot\nabla\psi,w_i)\in\real\,,\ i=1,2.
$$
Then, there is a positive constant $c$ depending, at most, on $T$, such that 
\be\ba{ll}\medskip
|\mathscr T_1|\le c\,\|\alpha\|_{L^{r,s}}\|w_1\|_{L^{r,s}}\|\psi\|_{\mathscr W^{r',s'}}\,, & \mbox{if $\alpha,w_1$ satisfy {\rm (i)} and $\psi\in \mathscr W^{r',s'}$}\,,\\
|\mathscr T_1|\le c\,\|\alpha\|_{L^{\infty,3}}\|w_1\|_{L^{\infty,3}}\|\psi\|_{\mathscr W^{3,\frac32}}\,, & \mbox{if $\alpha,w_1$ satisfy {\rm (ii)} and $\psi\in \mathscr W^{\frac32}$}\,;
\ea
\eeq{11}
and
\be\ba{ll}\medskip
|\mathscr T_2|\le c\,\|\alpha\|_{L^{r,s}}\|w_2\|_{L^{\infty,2}}\|\psi\|_{\mathscr W^{2}}\,, & \mbox{if $\alpha$ satisfies {\rm (i)} and $\psi\in \mathscr W^{2}$}\,,\\
|\mathscr T_2|\le c\,\|\alpha\|_{L^{\infty,3}}\|w_2\|_{L^{\infty,2}}\|\psi\|_{\mathscr W^{ 2}}\,, & \mbox{if $\alpha $ satisfy {\rm (ii)} and $\psi\in \mathscr W^{2}$}\,.
\ea
\eeq{110}
\EL{RoNa}
{\em Proof.} By the H\"older inequality,
$$
|\mathscr T_1|\le \|\alpha\|_{L^{r,s}}\|w_1\|_{L^{r,s}}\|\nabla\psi\|_{L^{\frac{r}{r-2},\frac s{s-2}}}\,.
$$
Since $2(r-2)/r+3(s-2)/s=2/r'+3/s'-1$, in the previous inequality  we may use Lemma~A.1 in the Appendix with $p\equiv r'$, $q\equiv s'$, which entails \eqref{11}$_1$. Moreover, by a further use of H\"older inequality, 
$$
|\mathscr T_1|\le \|\alpha\|_{L^{\infty,3}}\|w_1\|_{L^{\infty,3}}\|\nabla\psi\|_{L^{1,3}}\,, 
$$
so that \eqref{11}$_2$ follows from the latter and the general Sobolev inequality 
\be
\|\nabla\psi\|_{\frac{3p}{3-p}}\le c\,\|D^2\psi\|_{p}\,,\ \ p\in [1,3) 
\eeq{si}
with the choice $p=3/2$.
Likewise,
$$
|\mathscr T_2|\le \|\alpha\|_{L^{\infty,3}}\|w_2\|_{L^{\infty,2}}\|\nabla\psi\|_{L^{1,6}}\,, 
$$
and \eqref{110}$_2$ follows from the latter and \eqref{si} with $p=2$. It remains to show \eqref{110}$_1$. To this end, we observe that, by the H\"older inequality,
\be
|\mathscr T_2|\le \|\alpha\|_{L^{r,s}}\|w_2\|_{L^{\infty,2}}\|\nabla\psi\|_{L^{r',\frac{2s'}{2-s'}}}\,.
\eeq{ch}
Since $s'\in (1,3/2)$, it follows $2<2s'/(2-s')<6$, and so,
by interpolation,
$$
\|\nabla\psi\|_{\frac{2s'}{2-s'}}\le \|\nabla\psi\|_2^\theta\,\|\nabla\psi\|_6^{1-\theta}\,,\ \theta=\theta(s')\in(0,1)\,.
$$
Employing the continuous embedding $\mathscr W^{2}\subset C([0,T];H^1)$ \cite{Lions1}, the desired property is then a consequence of the latter inequality, \eqref{si} (with $p=2$), and \eqref{ch}.
\par\hfill$\square$
\Bl Let $v\in L^2_{{\rm loc}}(\real^3\times [0,T))$. Then, for any $v_0\in L^2_\sigma(\real^3)$ there exists $u\in C_w([0,T];L^2_\sigma(\real^3))\cap L^2(0,T;H^1(\real^3))$ such that
\be
\Int{0}{T}\left(u,\partial_t\varphi+\Delta\varphi+v\cdot\nabla\varphi\right)=-(v_0\cdot\varphi(0))\,,
\eeq{ABI_1}
for all $\varphi\in\mathcal D_T$. Moreover $\lim_{t\to0^+}\|u(t)-v_0\|_2=0$, and $u$ satisfies \eqref{1.3} (with $u\equiv v$).
\EL{ABI}
{\em Proof.} The result is achieved in an entirely routine fashion, by employing the classical Galerkin method; see, e.g., the proof given in \cite[Lemma 5.4]{Ga}.\par\hfill$\square$

\Bl Let  $\alpha$  satisfy  one of the following assumptions 
$$
{\rm (i)}\,\ \alpha \in L^{r,s}\,,\ r:=\frac{2s}{s-3}\,,\ s\in(3,\infty)\,;\ \ {\rm (ii)}\,\ \alpha \in C([0,T];L^3)\,,
$$
and let $F\in C_0^{\infty}(\real^3\times (0,T))$. Then, the problem
\be\ba{cc}\smallskip\left.\ba{ll}\smallskip
\partial_t\Psi+\alpha\cdot\nabla \Psi=\Delta \Psi-\nabla {\chi}+F\\
\Div \Psi=0\ea\right\}\ \ \mbox{in}\ \,\real^3\times (0,T)
\\
\Psi(x,0)=0\ \ x\in\real^3\,,
\ea
\eeq{2.1_0}
has one (and only one) solution $(\Psi,\chi)$ with
\be 
\Psi\in \mathscr W^{r',s'}\cap \mathscr W^2\,,\ \ \nabla\chi\in L^{r',s'}\cap L^{2,2}\,,\ \ \mbox{in case {\rm (i)}}\,,
\eeq{A2}
and
\be 
\Psi\in \mathscr W^{\frac32}\cap \mathscr W^2\,,\ \ \nabla\chi\in L^{\frac32,\frac32}\cap L^{2,2}\,,\ \ \mbox{in case {\rm (ii)}}\,.
\eeq{A3}
\EL{2.4}
{\em Proof.}  
We begin to consider the following regularized version of \eqref{2.1_0}
\be\ba{cc}\smallskip\left.\ba{ll}\smallskip
\partial_t\Psi+\alpha_{(\eta)}\cdot\nabla \Psi=\Delta \Psi-\nabla {\chi}+F\\
\Div \Psi=0\ea\right\}\ \ \mbox{in}\ \,\real^3\times (0,T)
\\
\Psi(x,0)=0\ \ x\in\real^3\,.
\ea
\eeq{2.1_R}
Since $C_0^\infty(\real^3\times(0,T))$ is dense in $L^{r',s'}$, for any given $\varepsilon>0$ we may write $\alpha=\alpha_1+\alpha_2$ where
\be
\|\alpha_1\|_{L^{r,s}}<\varepsilon\,,\ \ \|\alpha_2\|_{L^{\infty,\infty}}\le c\,\|\alpha\|_{L^{r,s}}\,.
\eeq{2.25}
Moreover, for $\lambda>0$ let 
\be\zeta={\rm e}^{-\lambda\,t}\Psi\,, \ \rho={\rm e}^{-\lambda\,t}\chi\,,\ G={\rm e}^{-\lambda\,t}F\,.
\eeq{2.26} 
As a consequence, problem \eqref{2.1_R} becomes
\be\ba{cc}\smallskip\left.\ba{ll}\smallskip
\partial_t\zeta+\lambda\,\zeta +({ {\alpha_1}_{(\eta)}}+{\alpha_2}_{(\eta)})\cdot\nabla\zeta=\Delta \zeta-\nabla {\rho}+G\\
\Div \zeta=0\ea\right\}\ \ \mbox{in}\, \real^3\times (0,T)
\\
\zeta(x,0)=0\ \ x\in\real^3\,,
\ea
\eeq{2.27} 
The existence of a solution $\zeta$ in the class $\mathscr W^{r',s'}$ can be established for $\lambda$ sufficiently large, by a simple perturbation argument around the solution to the problem obtained by formally setting $\alpha_1\equiv \alpha_2\equiv 0$ in \eqref{2.27}.  To show this, let
\be
\mathcal F:=({ {\alpha_1}_{(\eta)}}+{{\alpha_2}_{(\eta)}})\cdot\nabla\bar{\zeta}+G\,,\ \ \bar{\zeta}\in\mathscr W^{r',s'}\,,
\eeq{2.28}
so that \eqref{2.27} can be written as
\be\ba{cc}\smallskip\left.\ba{ll}\smallskip
\partial_t\zeta+\lambda\,\zeta =\Delta \zeta-\nabla {\rho}+\mathcal F\\
\Div \zeta=0\ea\right\}\ \ \mbox{in}\, \real^3\times (0,T)
\\
\zeta(x,0)=0\ \ x\in\real^3\,,
\ea
\eeq{2.29}
By well known results  (e.g. \cite[Theorem 5.5]{DHP}) one infers that \eqref{2.29} has one and only one solution $\zeta\in \mathscr W^{r',s'}$, $\nabla\rho\in L^{r',s'}$ such that  
$$
\int_0^T\big(\|\partial_\tau \zeta(\tau)\|_{s'}^{r'}+\lambda^{r'}\|\zeta(\tau)\|_{s'}^{r'}+\|D^2\zeta(\tau)\|_{s'}^{r'}\big)\le c\int_0^T\|\mathcal F(\tau)\|_{s'}^{r'}\,,
$$
with $c$ {\em independent} of $\lambda$. Thus, by choosing $\lambda\ge 2$ (say) along with classical interpolation inequalities, from the last displayed equation we deduce
\be
\|\zeta\|_{\mathscr W^{r',s'}}+\lambda\,\|\zeta\|_{L^{r',s'}}\le c\,\|\mathcal F\|_{L^{r',s'}}\,.
\eeq{2.30}
Next, by H\"older inequality, 
$$
\|{ {\alpha_1}_{(\eta)}}\cdot\nabla \bar{\zeta}\|_{L^{r',s'}}\le \|{ {\alpha_1}_{(\eta)}}\|_{L^{r,s}}\|\nabla \bar{\zeta}\|_{L^{\frac r{r-2},\frac s{s-2}}}\le \|\alpha_1\|_{L^{r,s}}\|\nabla \bar{\zeta}\|_{L^{\frac r{r-2},\frac s{s-2}}}\,.
$$
Since $2(r-2)/r+3(s-2)/s=2/r'+3/s'-1$, we may use in the previous inequality the embedding Lemma A.1 with $p\equiv r'$, $q\equiv s'$ along with  \eqref{2.25} to show
\be
\|{ {\alpha_1}_{(\eta)}}\cdot\nabla \bar{\zeta}\|_{L^{r',s'}}\le c\, \varepsilon \,\|\bar{\zeta}\|_{\mathscr W^{r',s'}}\,.
\eeq{2.31}
Furthermore, again by \eqref{2.25}, we infer
$$
\|{{\alpha_2}_{(\eta)}}\cdot\nabla\bar{\zeta}\|_{L^{r',s'}}\le \|{{\alpha_2}_{(\eta)}}\|_{L^{\infty,\infty}}\|\nabla\bar{\zeta}\|_{L^{r',s'}}\le c\,\|\alpha\|_{L^{r,s}}\|\nabla\bar{\zeta}\|_{L^{r',s'}}\,, 
$$
and so, using in this relation the Ehrling inequality
$$
\|\nabla\bar{\zeta}\|_{s'}\le\varepsilon\, \|D^2\bar{\zeta}\|_{s'}+c_\varepsilon\,\|\bar{\zeta}\|_{s'}
$$
we conclude
\be
\|{{\alpha_2}_{(\eta)}}\cdot\nabla\bar{\zeta}\|_{L^{r',s'}}\le c\,\|\alpha\|_{L^{r,s}}\big[\,\varepsilon\,\|\bar{\zeta}\|_{\mathscr W^{r',s'}}+\big(\frac{c_\varepsilon}\lambda\big) \lambda\,\|\bar{\zeta}\|_{L^{r',s'}}\big]\,.
\eeq{2.32}
Consider now the map $M:\bar{\zeta}\in W^{r',s'}\mapsto {\zeta}\in W^{r',s'}$ with $\zeta$ solving \eqref{2.29}, and endow $W^{r',s'}$ with the (equivalent) norm $\|\cdot\|_{\mathscr W^{r',s'}}+\lambda \,\|\cdot\|_{L^{r',s'}}$. From  
\eqref{2.28}, \eqref{2.30}--\eqref{2.32}, by taking $\varepsilon$ sufficiently small and $\lambda$ sufficiently large compared to $\|\alpha\|_{L^{r,s}}$, we  show at once that $M$ possesses a fixed point that solves \eqref{2.27}, and, in addition,  satisfies the estimate 
$$
\|\zeta\|_{\mathscr W^{r',s'}}+\lambda\,\|\zeta\|_{L^{r',s'}}\le c\,\|G\|_{L^{r',s'}}\,,
$$
with $c$ independent of $\eta$.
As a result, in view of \eqref{2.26} we obtain that problem  \eqref{2.1_R} has a solution $(\Psi_\eta,\chi_\eta)$ with $\Psi_\eta\in \mathscr W^{r',s'}$, $\nabla\chi_\eta\in L^{r',s'}$ that obeys the following estimate, uniformly in in $\eta$,
\be
\|\Psi_\eta\|_{\mathscr W^{r',s'}}\le c\,\|F\|_{L^{r',s'}}\,.
\eeq{2.33}
We now show by a simple boot-strap argument that, in fact, $\Psi_\eta\in \mathscr W^{2}$. Since, by Lemma A.1, $\nabla\Psi_\eta\in L^{\frac r{r-2},\frac s{s-2}}$ and ${\alpha_{(\eta)}}\in L^{\infty,\infty}$ we deduce by classical existence and uniqueness theory for problem \eqref{2.1_R} (e.g. \cite[Theorem VIII.4.1 and Lemma VIII.4.2]{GaBook}) that $\Psi_\eta\in \mathscr W^{\frac r{r-2},\frac s{s-2}}$. By Lemma A.1, this property entails $\nabla \Psi_\eta\in L^{2,\frac43}$ which, in turn, ensures $\Psi_\eta\in \mathscr W^{2,\frac43}$. Again by Lemma A.1, the latter furnishes $\nabla\Psi_\eta\in L^{4,2}$ which, finally, implies $\Psi_\eta\in \mathscr W^{2}$. With this information in hand, it is then a routine task (see, e.g., \cite[Lemma 5.4]{Ga}) to show that $\Psi_\eta$ is bounded in $\mathscr W^{2}$, uniformly in $\eta$. We will thus only sketch the proof. By testing \eqref{2.1_R}, in the order,  with $\Psi_\eta$, $\Delta\Psi_\eta$, and $\partial_t\Psi_\eta$, we get
\be\ba{ll}\medskip
\half \ode{}t\|\Psi_\eta\|_2^2+\|\nabla\Psi_\eta\|_2^2=\big(F,\Psi_\eta\big)\\ \medskip
\half \ode{}t\|\nabla\Psi_\eta\|_2^2+\|\Delta\Psi_\eta\|_2^2=\left(\alpha_{(\eta)}\cdot\nabla\Psi_\eta-F,\Delta\Psi\right)\\\
\half \ode{}t\|\nabla\Psi_\eta\|_2^2+\|\partial_t\Psi_\eta\|_2^2=-\left(\alpha_{(\eta)}\cdot\nabla\Psi_\eta-F,\partial_t\Psi_\eta\right)\,.
\ea
\eeq{abi}
By H\"older and Sobolev inequalities and the property of mollifiers, for any $\varepsilon>0$ one can show (see \cite[eq. (5.10)]{Ga})
\be
(\alpha_{(\eta)}\cdot\nabla\Psi_\eta,V)\le c_\varepsilon\,\|\alpha\|_{s}^r\|\nabla\Psi_\eta\|_2^2+\varepsilon\left(\|\Delta \Psi_\eta\|_2^2+\|V\|_2^2\right))\,,\ \ V:=\Delta\Psi_\eta\,,\ -\partial_t\Psi_\eta\,.
\eeq{gail}
Therefore, combining \eqref{abi}, \eqref{gail} and using also Cauchy-Schwarz inequality, with the help of Gronwall's lemma we conclude 
\be
\|\Psi_\eta\|_{\mathscr W^{2}}\le c\,\|F\|_2\,,
\eeq{estim}
with $c$ independent of $\eta$. The desired existence result in case (i) then follows by letting $\eta\to0$ (along a sequence) in \eqref{2.1_R},  and using the uniform estimates \eqref{2.33} and \eqref{estim}.
In the case (ii), we again start from the modified system \eqref{2.27} where, this time, $\alpha_1$ and $\alpha_2$ are chosen with the property
\be
\|\alpha_1\|_{L^{\infty,3}}<\varepsilon\,,\ \ \ \|\alpha_2\|_{L^{\infty,\infty}}\le c\,\|\alpha\|_{L^{\infty,3}}\,.
\eeq{2.35}
Proceeding as in the proof of case (i), we show that problem \eqref{2.29} has a solution $\zeta\in \mathscr W^{\frac32}$ such that
\be
\|\zeta\|_{\mathscr W^{\frac32}}+\lambda\,\|\zeta\|_{L^{\frac32,\frac32}}\le c\,\|\mathcal F\|_{L^{\frac32,\frac32}}\,.
\eeq{2.36}
with $\mathcal F$ as in \eqref{2.28}.
Now, for given $\bar{\zeta}\in \mathscr W^{{\frac32}}$, by H\"older inequality we get
$$
\|{ \alpha_1}_{(\eta)}\cdot\nabla \bar{\zeta}\|_{L^{\frac32,\frac32}}\le \|{ \alpha_1}_{(\eta)}\|_{L^{\infty,3}}\|\nabla \bar{\zeta}\|_{L^{\frac32,3}}\le \|\alpha_1\|_{L^{\infty,3}}\|\nabla \bar{\zeta}\|_{L^{\frac32,3}}.
$$
As a result, since  $\|\nabla\bar\zeta\|_3\le c\,\|D^2\bar\zeta\|_{\frac32}$ (see \eqref{si} with $p=3/2$), thanks to \eqref{2.35} we may deduce 
\be
\|{ {\alpha_1}_{(\eta)}}\cdot\nabla \bar{\zeta}\|_{L^{\frac32,\frac32}}\le c\,\varepsilon\,\|\bar\zeta\|_{\mathscr W^{\frac32}}\,.
\eeq{2.37}
Likewise, by an argument entirely analogous to that leading to \eqref{2.32}, which now uses \eqref{2.35}, we show
\be
\|{{\alpha_2}_{(\eta)}}\cdot\nabla\bar{\zeta}\|_{L^{\frac32,\frac32}}\le c\,\|\alpha\|_{L^{\infty,3}}\big[\,\varepsilon\,\|\bar{\zeta}\|_{\mathscr W^{\frac32}}+\big(\frac{c_\varepsilon}\lambda\big) \lambda\,\|\bar{\zeta}\|_{L^{\frac32,\frac32}}\big]\,.
\eeq{2.38}
Employing \eqref{2.36}--\eqref{2.38} and a fixed-point procedure as in the proof of case (i), we conclude that \eqref{2.27} has a solution $(\Psi_\eta,\chi_\eta)$  with
$ 
\Psi_\eta\in \mathscr W^{\frac32}$, $\nabla\chi_\eta\in L^{\frac32,\frac32}$ that obeys the following estimate, uniformly in $\eta$
\be
\|\Psi_\eta\|_{\mathscr W^{\frac32}}\le c\,\|F\|_{L^{\frac32,\frac32}}\,.
\eeq{gb}
Next, from Lemma A.1 with $p=q=3/2$ it follows $\nabla\Psi_\eta\in L^{\frac{12}5,2}$, and so, since $\alpha_{(\eta)}\in L^{\infty,\infty}$, by existence and uniqueness theory for problem \eqref{2.1_R} (e.g. \cite[Theorem VIII.4.1 and Lemma VIII.4.2]{GaBook}) we get $\Psi_\eta\in \mathscr W^{2}$. Once this property has been established, we can proceed exactly as in the proof of case (i) to obtain the uniform bound \eqref{estim}, also for the case at hand. To this end, it suffices to replace \eqref{gail} with the following one (see \cite[Lemma 5.3]{Ga})
$$
(\alpha_{(\eta)}\cdot\nabla\Psi_\eta,V)\le c_\varepsilon\,\|\alpha\|_3^2\|\nabla\Psi_\eta\|_2^2+\varepsilon\,\left(\|\Delta\Psi_\eta\|_2^2+\|V\|_2^2\right)\,,\ \ V:=\Delta\Psi_\eta\,,\ -\partial_t\Psi_\eta\,.
$$ 
The proof then is achieved by letting $\eta\to0$ (along a sequence) in \eqref{2.1_R} and using the uniform estimates \eqref{gb} and \eqref{estim}.
\par\hfill$\square$

\setcounter{equation}{0}
\section{Proofs of Theorem 1.1 and its Corollaries}
{\em Proof of \theoref{1}}. We shall only detail the proof in the case \eqref{H}$_1$, since the case \eqref{H}$_2$ is treated in an entirely similar way. 
Let $\chi_\beta=\chi_\beta(t)$, $3\beta\in (0,T)$ be a smooth non-negative function of $t$ such that $\chi_\beta(t)=1$ if $t\ge3\beta$, and $=0$, if $t\in [0,2\beta]$, with $|\chi'_\beta(t)|\le c/\beta$. Replacing $\varphi$ with $\chi_\delta\varphi$ in both \eqref{1.5} and \eqref{ABI_1}, and setting $w=v-u$, with $u$ given in \lemmref{ABI}, we obtain
for arbitrary $\varphi\in \mathcal D_T$,
\be
\Int{2\delta}{T}\big(\chi_\delta\,w,\partial_t\varphi+\Delta\varphi+\chi_{\frac{2\delta}3}v\cdot\nabla\varphi\big)=-\int_{2\delta}^{3\delta}\chi'_\delta(w,\varphi)\,,
\eeq{4.1}
where we also have used $\chi_\delta(t)\chi_{\frac{2\delta}3}(t)=\chi_\delta(t)$, $t\in [0,T]$.
Thus,  $\chi_\delta\,w\in L^{r,s}+L^{\infty,2}$ and $\chi_{\frac{2\delta}3}\,v\in L^{r,s}$. As a result, from  \lemmref{N1} and \lemmref{RoNa} 
we can readily show that \eqref{4.1} leads to the following one
\be
\Int{2\delta}{T}\big(\chi_\delta\,w,\partial_t\psi+\Delta\psi+\chi_{\frac{2\delta}3}v\cdot\nabla\psi\big)=-\int_{2\delta}^{3\delta}\chi'_\delta(w,\psi)\,,\ \ \mbox{for all $\psi\in \mathscr W^{r',s'}_0\cap\mathscr W^2_0$}\,.
\eeq{4.1_0}
Let  $\psi(x,t):=\Psi(x,T-t)$,   where $\Psi$  is the solution to \eqref{2.1_0} constructed in \lemmref{2.4} with $\alpha(x,t):=-(\chi_{\frac{2\delta}3}\,v)(x,T-t)$ and $F(x,t):=-f(x,T-t)$, $(x,t)\in\real\times [0,T]$, with $f\in C_0^\infty(\real^3\times(0,T))$. By that lemma it then follows that $\psi$ is in the class $\mathscr W^{r',s'}_0\cap\mathscr W^2_0$, and can be thus used as test function in \eqref{4.1_0}. Therefore,  also with the help of \eqref{2.1_0}$_1$, we conclude
\be
\Int{2\delta}{T}\big(\chi_\delta\,w,f\big)=-\int_{2\delta}^{3\delta}\chi'_\delta(w,\psi)\,.
\eeq{4.1_00}
Clearly,
\be
\lim_{\delta\to 0}\Int{2\delta}{T}\big(\chi_\delta\,w,f\big)=\Int{0}{T}\big(w,f\big)\,.
\eeq{4.1_000}
Furthermore, since by assumption and \lemmref{ABI}, $w(t)\to 0$ weakly in $L^2_\sigma(\real^3)$ and, by \lemmref{2.4} and  the continuous embedding $\mathscr W^2\subset C([0,T];H^1)$, $\psi(t)\to \psi(T)$ strongly in $H^1(\real^3)$, recalling that $\chi_\delta^\prime\sim \delta^{-1}$, we also have
$$
\lim_{\delta\to0}\int_{2\delta}^{3\delta}\chi'_\delta(w,\psi)=0\,.  
$$
From the latter, \eqref{4.1_00} and \eqref{4.1_000} we infer
$$
\Int{0}{T}\big(w,f\big)=0\,,\ \ \mbox{for all $f\in C_0^\infty(\real^3\times(0,T))$}\,,
$$which implies $w\equiv 0$, namely, $v\equiv u$.
The  theorem is therefore completely proved. \par\hfill$\square$\par 
{\em Proof of \cororef{1.1}}. From \theoref{1} and classical results (e.g. \cite[Theorem 5]{Serrin}) we show, for all sufficiently small $\delta>0$, 
\be
\|v(t)\|_2^2+2\int_\delta^t\|\nabla v(\tau)\|_2^2=\|v(\delta)\|_2^2\,,\ \ \mbox{all $t\ge \delta$}\,,
\eeq{PI}
and the result follows by letting $\delta\to 0$ in the latter and, again, using \theoref{1}.\par\hfill$\square$
\par 
{\em Proof of \cororef{1.1_0}}. It is an immediate consequence of \cororef{1.1}.\par\hfill$\square$
\par 
{\em Proof of \cororef{U}}. It follows from \theoref{1} and classical results concerning the uniqueness of Leray-Hopf solutions; see, e.g. \cite[Theorem 4.2 and Theorem 7.2]{Ga}.\par\hfill$\square$
\par 
{\em Proof of \cororef{1.2}}. Under the assumption \eqref{H}$_1$, the result follows from \theoref{1} and \cite[Theorem 2]{MM}. We will give a straightforward proof that covers both conditions in \eqref{H}. By \theoref{1} we have, in particular, that $v$ is a Leray-Hopf solution on $[\delta,T]$,  which, by assumption satisfies one of the properties in \eqref{H}. This then implies  $v\in C^\infty(\real^3\times(\delta,T])$ (see, e.g., \cite[Theorem 5.2]{Ga}). Since $\delta\,(>0)$ is arbitrary, the result is proved.\par\hfill$\square$
\Br It is likely that \theoref{1} may continue to hold in the borderline situation $v\in L^{\infty,3}$. However, it is not clear whether the method presented here would work in that case.
\ER{1}
\Br It is easy to check that the proof of the key \lemmref{2.4} and, therefore, of \theoref{1}, continues to be valid also in space dimension $n=4$. However, if $n\ge 5$ some of the embedding inequalities used in that lemma no longer hold, thus leaving the validity of \theoref{1} open in such a case.  
\ER{2}
\section*{Appendix}
\renewcommand{\theequation}{A.\arabic{equation}}\setcounter{equation}{0}\noindent{\bf Lemma A.1} {\sl If $u\in \mathscr W^{p,q}$, then
\be
\nabla u\in L^{p_1,q_1}\,,\ \ \frac2{p_1}+\frac3{q_1}=\frac2{p}+\frac3{q}-1\,.
\eeq{E2}
}\par
{\em Proof.} The proof  can be obtained by an argument similar to that used in \cite[Theorem 2.1]{Solo}. More precisely, we start with the classical representation of  $u\in \mathscr W^{p,q}$
\be\ba{rl}\medskip
u(x,t)=&\Int{t-1}t\left(\Int{\real^2}{}\Gamma_1(x-y,t-\tau)(u_\tau-\Delta u)(y,\tau)dy\right)d\tau\\&+\Int{t-1}t\left(\Int{\real^2}{}\Gamma_2(x-y,t-\tau)u(y,\tau)dy\right)d\tau\,,
\ea
\eeq{p_7}
were $\Gamma_1(x,t)=\Gamma (x,t)\psi(|x|)\psi(t)$, while
$$
\Gamma_2(x,t)=2\psi(t)\nabla\Gamma(x,t)\cdot\nabla\psi(|x|)+\Gamma(x,t)\left(\psi(t)\Delta\psi(|x|)-\psi'(t)\psi(|x|)\right)\,,
$$
with 
\be\Gamma(x,t)=\left\{\ba{ll}\smallskip
\Frac{1}{4\pi t}\exp\left(-\Frac{|x|^2}{4t}\right)&\ \ \mbox{if $t>0$}\\
0 &\ \ \mbox{if $t<0$}\ea\right.\,,
\eeq{p_8}
and $\psi:\real\to\real$ a smooth non-negative function such that $\psi\le 1$, $\psi(\xi)=0$ if $\xi\ge 1$ and $\psi(\xi)=1$ if $\xi\le 1/2$. By using  in \eqref{p_7} the generalized Minkowski inequality followed by Young's inequality for convolutions we get 
$$
\|\nabla u(t)\|_{q_1}\!\le\!\int_{t-1}^t\!\!\!\big(\|\nabla\Gamma_1(t-\tau)\|_{\sigma}\!+\|\nabla\Gamma_2(t-\tau)\|_{\sigma}\big)\big(\|u_\tau(\tau)\|_{q}\!+\|u(\tau)\|_{2,q}\big)d\tau,
$$ 
where $q_1\ge q$, and $1/\sigma=1+1/q_1-1/q$. As a result, taking into account that (as shown by an easy calculation that uses \eqref{p_8} and the properties of the function $\psi$)
$$
|\nabla\Gamma_i(\xi,s)|\le \Frac{c_1}{(|\xi|^2+s)^{2}}\,,\ \ \ i=1,2\,,
$$
we find
\be
\|\nabla u(t)\|_{q_1}\le c_2\int_{t-1}^t\Frac{\|u_\tau(\tau)\|_{q}+\|u(\tau)\|_{2,q}}{(t-\tau)^{\frac32(\frac1q-\frac1{q_1})+\frac12}}d\tau\,,
\eeq{p_9}
with $c_2$ independent of $t$. The property stated in \eqref{E2} then follows by employing in \eqref{p_9}  the Hardy-Littlewood inequality.\par\hfill$\square$

\ed